\documentclass[10pt,conference]{IEEEtran} 
\IEEEoverridecommandlockouts
\usepackage{cite}
\usepackage{amsmath,amssymb,amsfonts}
\usepackage{algorithmic}
\usepackage{graphicx}
\usepackage{textcomp}
\usepackage{xcolor}
\usepackage{url}
\usepackage{commath} 

\newtheorem{thm}{Theorem}

\def\BibTeX{{\rm B\kern-.05em{\sc i\kern-.025em b}\kern-.08em
    T\kern-.1667em\lower.7ex\hbox{E}\kern-.125emX}}
\begin{document}

\title{The L\'{e}vy State Space Model}

\author{\IEEEauthorblockN{Simon~Godsill}
\IEEEauthorblockA{\textit{Department of Engineering} \\
\textit{University of Cambridge}\\
Cambridge, CB2 1PZ, UK \\
sjg30@cam.ac.uk}
\and
\IEEEauthorblockN{Marina Riabiz}
\IEEEauthorblockA{
\textit{Department of Biomedical Engineering} 
\\
\textit{King's College London}\\
London, SE1 7EH, UK \\
marina.riabiz@kcl.ac.uk}
\and
\IEEEauthorblockN{Ioannis Kontoyiannis}
\IEEEauthorblockA{\textit{Department of Engineering} \\
	\textit{University of Cambridge}\\
	Cambridge, CB2 1PZ, UK \\
	ik355@cam.ac.uk}
}

\maketitle


\begin{abstract}
In this paper we introduce a new class of state space models based on shot-noise simulation representations of non-Gaussian L\'{e}vy-driven linear systems, represented as stochastic differential equations. In particular a conditionally Gaussian version of the models is proposed that is able to capture heavy-tailed non-Gaussianity while retaining tractability for inference procedures. We focus on a canonical class of such processes, the $\alpha$-stable L\'{e}vy processes, which retain important properties such as self-similarity and heavy-tails, while emphasizing that broader classes of non-Gaussian L\'{e}vy processes may be handled by similar
methodology. An important feature is that we are able to marginalise both the skewness and the scale parameters of these challenging models from posterior probability distributions. The models are posed in continuous time and so are able to deal with irregular data arrival times. Example modelling and inference procedures are provided using Rao-Blackwellised sequential Monte Carlo applied to a two-dimensional Langevin model, and this is tested on real exchange rate data.
\end{abstract}


\section{Introduction}

Heavy-tailed non-Gaussian processes play a significant role in many real-world scenarios, exhibiting extreme values 
much more frequently than a Gaussian model. 
Examples of such abrupt changes include variations presented by stock prices 
or insurance gains/losses in financial applications, 
as studied extensively since the seminal works 
\cite{Mandelbrot1963} 
and \cite{Fama1965a}.
Further applications can be found in various fields of engineering, such as communications (see \cite{Azzaoui2010} for statistical modelling of channels,  \cite{Fahs2012, FreitasEganClavierEtAl2017} for capacity bounds, 
\cite{LiebeherrBurchardCiucu2012} for delay bounds in networks 
with $\alpha$-stable noise, and \cite{ShevlyakovKim2006, WarrenThomas1991} 
for signal detection), 
signal processing \cite{Nikias1995,Lombardi2006,GodsillRayner1998,GodsillKuruoglu1999,Godsill_2000a}, image analysis \cite{Achim2001, Achim2006}.

In many of these situations, the random phenomena considered can be still thought 
of as emerging from the combination of many independent perturbations. 
According to the generalized CLT 
\cite[p.~162]{Gnedenko1968}\cite[p.~576]{Feller1966}, 
whenever the sum of independent identically distributed (i.i.d.) random 
variables (RVs) converges in distribution, it converges to 
a member of the class of $\alpha$-stable distributions;
this class is central to this paper. The Gaussian 
is a special member of this class, the only one with finite variance. 
Hence, using non-Gaussian $\alpha$-stable processes offers a natural way to extend beyond Gaussian processes to heavy-tailed cases.
A principal reason for the attention that $\alpha$-stable laws have received
in applications (see the extensive bibliography listed 
in \cite{nolan_web}),
is their role in the generalized CLT, the modelling flexibility offered by this  class of distributions, as well as the useful and unique properties of $\alpha$-stable L\'{e}vy processes \cite{Samorodnitsky_Taqqu_1994} such as self-similarity and infinite divisibility. We will thus focus on the $\alpha$-stable class of models, but will also describe how our models and methods can be applied to a much broader class of heavy-tailed processes. 


\section{Problem formulation}
We study in particular inference methods for linear stochastic differential equations (sdes) driven by non-Gaussian L\'{e}vy processes 
\cite{rogers_williams_2000,Oksendal2003,Bertoin_1997}: 
\begin{IEEEeqnarray}{lCr}
d{X}(t)=A{X}(t)dt+\mathbf{h}dW(u),\quad {X}(t)\in \Re^P,\,\,W(t)\in\Re 
\IEEEeqnarraynumspace
\label{eq:sde}
\end{IEEEeqnarray}
where
$\{W(t)\}$ is  a non-Gaussian L\'{e}vy process, whose solution is obtained by stochastic integration:
\[
{X}(t)=e^{At}{X}(0)+\int_{0}^{t}e^{A(t-u)}\mathbf{h}dW(u)\,.
\]
The characteristic function for a general  L\'{e}vy process having no drift or Brownian motion part is given by \cite{Kallenberg_2002}, Corollary 15.8, \cite{Bertoin_1997}, as
\begin{equation}
E\left[\exp(iuW(t)) \right]= \exp( t\psi(u))\end{equation}
where
\begin{equation} \psi(u)=\int_{\Re^M \backslash \{0\}}
(e^{iuw}
-1 - iuw{\cal I}(|w|\leq 1))Q(dw)
\end{equation}
where ${\cal I}(\cdot)$ is the indicator function
and $Q$ is a L\'{e}vy measure on $\Re^M{\backslash}\{0\}$. We will be particularly concerned with the so-called {\em infinite activity\/} processes having $\int Q(dw)=\infty$, and hence an almost surely infinite number of jumps within any finite time interval.

We will consider principally the $\alpha$-stable L\'{e}vy process, later highlighting how to broaden the analysis to other classes. The focus will be on vector sdes driven by scalar L\'{e}vy processes,  although in principle our methods extend also to vector valued L\'{e}vy processes.
Take, then, $\{W(t)\}$ to be an $\alpha$-stable 
L\'{e}vy process, defined to have, for $0<\alpha<2$  \cite{Samorodnitsky_Taqqu_1994},
\begin{itemize}
	\item $W(0)=0$ a.s.
	\item Independent, stationary  $\alpha$-stable incrememts $W(t)-W(s)\overset{d}{=}W_{t-s}\sim S_{\alpha}((t-s)^{1/\alpha},\beta,0)$, $t>s$
\end{itemize}
where $\overset{d}{=}$ denotes having the same probability distribution, and  $S_{\alpha}(\sigma,\beta,\mu)$ is the $\alpha$-stable law, having scale parameter $\sigma$, skewness $\beta$, location $\mu$, and tail parameter $\alpha$ \cite{Feller1966,Samorodnitsky_Taqqu_1994,Nikias1995}.  
Such L\'{e}vy processes are pure jump processes ($0<\alpha<1$), or pure jump plus drift ($1\leq \alpha <2$), possessing 
almost surely an infinte number of jumps in finite time intervals and they are in general highly intractable for inference.\footnote{We explicity exclude the special case $\alpha=1$ for notational reasons, but this could in principle 
	be included in future studies.}

It will be convenient to define the integral, over the time interval $[0,1]$\footnote{We take the time axis up to $t=1$ by convention, although the processes can be extended to arbitrarily large times axes in practice.}
\[{\bf I}({\bf f}_t)=\int_{0}^{1}{\bf f}_tdW(u)
\]
with 
\[{\bf f}_t=e^{A(t-u)}{\cal I}(u\leq t)\mathbf{h}.\]

In the case $P=1$ (univariate sde)  with $A=a$, $\mathbf{h}=h$, the problem is solved using results based on \cite{Samorodnitsky_Taqqu_1994} Section 3.2.2, see our previous work \cite{Lemke2015,Lemke_Godsill_2015}.
In the vector sde case, $P>1$ we could consider for example a direct Euler-type discretisation of the integral, solving approximately using increments of the process, $\delta W_{j}= W((j+1)\delta t)-W(j\delta t)$, $j=0,1,...,N-1$, $\delta_t=t/N$,  and then using the conditionally Gaussian auxiliary variable models of \cite{Lemke2015,godsill2006bayesian}. However, there are issues about the level of approximation, choice of discretisation step-size, etc. that make this approach unsatisfactory as a general solution. 
However, a  pertinent question remains as to whether the limit as $N\rightarrow\infty$
can be explicitly calculated, hence avoiding the need to
simulate auxiliary states on the $\delta t$ scale. We have been unable to do so, except for the scalar $P=1$ case (see   \cite{Lemke_Godsill_2015}), and hence we propose very accurate approximate models suitable for likelihood and Bayesian inference procedures. 


\section{Summary of Methods and Contribution}
The starting point for the L\'{e}vy state space model is the generalised shot noise representation of $\{W(t)\}$, see \cite{Rosinski_2001}. A principal contribution here is to adopt a non-centered {\em conditionally Gaussian\/} form for the shot noise representation which has been largely overlooked in the stochastic processes literature, but which is of significant benefit in generating tractable and effective inference procedures.
This special form enables both tractable conditional Gaussian inference, and the possibility of modelling of asymmetric L\'{e}vy processes. In addition, we show that the residual error committed by truncation can be approximated by a Gaussian-driven sde that exactly matches the first and second moments of the true residual error. We have established central limit theorems and sharp rates of convergence for the approximating underlying Gaussian residual L\'{e}vy process. 
Finally, based on this representation, we develop transition densities for the {\em L\'{e}vy state space\/} model in the conditionally Gaussian form, and show how to perform marginal inference in the model using Kalman-filtering recursions built into a Bayesian (or potentially likelihood-based) Monte Carlo framework. Results have been obtained for simulated data and real exchange rate data.

Other relevant work includes \cite{clement2019} who propose theorems and estimating function methods for stable L\'{e}vy-driven sdes; \cite{Fournier_20011}, who prove error bounds between L\'{e}vy-driven sdes with small jumps and their Gaussian approximations; \cite{Jasra2019} provide Bayesian inference procedures using a quasi-likelihood approach and MCMC; \cite{Dereich2011,Dereich2016} provide a Multi-level Monte Carlo approach for evaluation of expectations, using coupled Euler approximations; and \cite{Ferreiro-Castilla2016} who propose Euler schemes under Poisson-arrival times. In contrast to these approaches, we tackle the problem without resort to Euler approximations or pseudo-likelihood functions, and our methods are applicable to low- or high-frequency observations, but are currently limited to linear sdes. Our approach focuses on a class of L\'{e}vy processes that may be expressed as a conditionally Gaussian shot noise process, which includes for example all $\alpha$-stable L\'{e}vy processes and many truncated and modified variants on these processes; however, at the expense of more computational effort, more general non-Gaussian versions may also be incorporated into our framework.   

\section{Generalised Shot Noise Processes}
We will use a generalised shot noise formulation studied in \cite{Rosinski_2001} and references therein, in which the L\'{e}vy process is expressed on a time axis $t\in[0,1]$ as
\begin{equation}
W(t)=\sum_{i=1}^\infty H(\Gamma_i,U_i){\cal I}(V_i\leq t)-tc_i.\label{eq:shot_noise}
\end{equation}
Here $\{\Gamma_i\}_{i=1}^\infty$ are the epochs of a unit rate Poisson process, $\{U_i\}_{i=1}^\infty$, $U_i\in\Re^P$ are mutually independent random variables with a certain distribution $F_U$, and $V_i$ are independent uniform random variables on $[0,1]$. The centering terms $c_i$ are sometimes needed and sometimes not, as discussed below.  $H(\Gamma,U)$ is generally taken to be {\em non-increasing\/} in its argument $\Gamma$ for fixed $U$. 
A point process for the expanded space $(\Gamma_i,U_i,V_i)\in \Re\times \Re^P\times [0,1]$ with mean measure $Lebesgue\times F_U\times Lebesgue$ may be constructed as \cite{Rosinski_2001}:
\[
M=\sum_{i=1}^\infty \delta_{V_i,\Gamma_i,U_i} 
\]
where the iid uniform variables $\{V_i\in[0,1]\}$ give the times of arrival of jumps and $\{H(\Gamma_i,U_i)\}$ give the size of the jumps and $\delta_{V,\Gamma,U}$ is the Dirac point mass located at  $V,\Gamma,U$.
This
leads to a direct expression for $X(t)$ 
\cite{Rosinski_2001} (by the L\'{e}vy-Ito integral representation of $W(t)$, expressed in terms of the augmented point process $M$):
\begin{align}
W(t) &=
\int_{
	|H(\gamma,u)|\leq 1}
H(\gamma,u) [(M([0, t], d\gamma, du) - tF_U(du)]\nonumber\\& +
\int_{
	|H(\gamma,u)|>1}
H(\gamma,u)M([0, t], d\gamma,du)\label{eq:Levy_direct_M}
\end{align}
which, on substitution for $M$, leads to the {\em shot noise\/} representation (\ref{eq:shot_noise}).
The almost sure convergence of this series to $\{W(t)\}$ is proven in \cite{Rosinski_2001}. 
In this paper we will apply, in particular, random Poisson truncations of the form $\{\Gamma_i\}_{i:\Gamma_i\leq c}$ which leads to approximations, neglecting small jumps, of the form   
\begin{IEEEeqnarray}{lCr}
W^c(t)
= & & 
\nonumber
\\
\,\,\,\int_{
	\gamma\in[0,c],|H(\gamma,u)|\leq 1}
H(\gamma,u) [(M([0, t], d\gamma, du) - tF_U(du)]
\nonumber
\\
+
\,\,\,\int_{
	\gamma\in[0,c],|H(\gamma,u)|>1}
H(\gamma,u)M([0, t], d\gamma,du)
\IEEEeqnarraynumspace
\label{eq:Levy_direct_M_c}
\end{IEEEeqnarray}
and the corresponding truncated shot noise series
\begin{equation}
W^c(t)=\sum_{i:\Gamma_i\in[0,c]} H(\Gamma_i,U_i){\cal I}(V_i\leq t)-tA(c).\label{eq:truncated_shot}
\end{equation}
The corresponding L\'{e}vy measure for $\{W(t)\}$ can then be shown to be \cite{Rosinski_2001}
\begin{equation}
Q(.)=\int_0^\infty \sigma_H(\gamma,.)d\gamma \label{eq:cond_Gauss_levy} 
\end{equation}
with
\[
A(c)=\int_0^c\int_{|x|<1}x\sigma_H(\gamma,dx)d\gamma,\,\,\,c>0
\] 
and where
\[
\sigma_H(\gamma,.)=P(H(\gamma,U)\in.),\,\,\,\gamma>0\,.
\]


\section{The conditionally Gaussian L\'{e}vy process}
One of the principal contributions of this paper is in methods for tractable inference about states and parameters of L\'{e}vy driven sdes.  The starting point for this is the generalised shot noise representation (\ref{eq:shot_noise}). We propose to use a special form of this in which we define:
\[
H(\Gamma,U)=h(\Gamma)U
\]
where $h(\cdot)\geq 0$ is a non-increasing function. 

If we furthermore assume a non-centered Gaussian form $U_i\sim {\cal N}(\mu_W,\sigma_W^2)$ then special properties result, notably tractable conditional Gaussian inference, and the possibility of modelling of asymmetric L\'{e}vy processes. Note that in the centered case, $\mu_W=0$, such processes are well known, see e.g. the processes of type G \cite{Rosinski_1991}, and the Condition S processes described in \cite{Samorodnitsky_Taqqu_1994}, but in the general case they are commonly overlooked. A well known example of this form for $H(\cdot)$ is the $\alpha$-stable L\'{e}vy process, which can be represented with the form \cite{Samorodnitsky_Taqqu_1994}:
\[
h(\Gamma)=\Gamma^{-1/\alpha}
\]
where $\alpha$ is the stable law tail parameter. 
The case $\alpha=1$, the Cauchy process, is specifically {\em excluded\/} from now on, as this forms a singular special case of the family \cite{Samorodnitsky_Taqqu_1994}. 
Essentially {\em any\/} distribution for the $U_i$ may be adopted, subject to finite absolute $\alpha$th  moments $E|U_i|^\alpha$, which is well known to be satisfied by our non-centred normal ${\cal N}(\mu_W,\sigma_W^2)$. Following the procedure in the previous section, it is relatively straightforward to verify that the correct L\'{e}vy measure for the stable law is obtained under this function $ h(\Gamma)=\Gamma^{-1/\alpha}$ and the non-centered Gaussian form $U_i\sim {\cal N}(\mu_W,\sigma_W^2)$ (a proof will be given in subsequent publications). 

There are many possible choices for the centering function and we follow a standard procedure \cite{Samorodnitsky_Taqqu_1994} in which no centering is required for $\alpha\in(0,1)$ and a modified centering can be calculated in closed form for $\alpha\in(1,2)$:
{\small
\[
A(c)=\int_{0}^c\int_{-\infty}^\infty xN(x|h(\gamma)\mu_W,h(\gamma)^2\sigma_W^2)dxd\gamma=\mu_W\frac{\alpha}{\alpha-1} c^{\frac{\alpha-1}{\alpha}}
\]}
and it can be shown that this choice of centering is also equal to the expected value of the truncated process:
\[
E\left[\sum_{i:\Gamma_i\in[0,c]} H(\Gamma_i,U_i){\cal I}(V_i\leq t)\right]=tA(c)=t\mu_W\frac{\alpha}{\alpha-1} c^{\frac{\alpha-1}{\alpha}}
\]
There is extensive earlier work on the analysis of the convergence 
rate of truncated series of this class to the corresponding stable laws; 
see, e.g., \cite{JanickiWeron1994,janicki1992computer,JanickiKokoszka1992,%
	LedouxPaulauskas1996, Bentkus1996, Bentkus2001}. 
In \cite{Riabiz2018a} we reviewed these
results, and derived new bounds for the symmetric stable laws under our conditionally Gaussian representation. 
Here though we are able to provide the following more general result that applies to the asymmetric general case for the first time:

\begin{thm}
	If the approximating truncated  process is $W^c(t)$ as in (\ref{eq:truncated_shot}) and $W(t)$ is the full untruncated L\'{e}vy process, then:
	\begin{IEEEeqnarray*}{L}
	\underset{x\in \Re}{\sup}|P(W^c(1)\leq x)-P(W(1)\leq x)|\leq
	\\
	\qquad\qquad\qquad     C(\sigma_W,\mu_W)\left(\frac{(2-\alpha)^{3/2}}{(3-\alpha)\alpha^{1/2}}\right)c^{-0.5}
	\end{IEEEeqnarray*}
	where $C(\sigma_W,\mu_W)$ is a constant that does not depend on the truncation level $c$.  
\end{thm}
{\em Proof:\/}
The proof follows a similar structure to Theorem 3.1 of \cite{Asmussen2001} and will be presented in a future publication.

{\em Remarks:\/}

\begin{enumerate}
	\item
	We notice that the rate of convergence wrt $c$ matches that in \cite{Asmussen2001} Th. 3.1, which corresponds to a different type of truncation on the absolute size of the jumps. We comment however that it is not as good as the convergence rate we have previously obtained for the {\em symmetric\/} case, which was $\sim c^{-1}$ \cite{Riabiz2018a}. However, we note that our current result is significantly more general in that it covers both symmetric and non-symmetric (skewed) cases. We can postulate that a more subtle future development will prove convergence rates between $c^{-0.5}$ and $c^{-1}$ as we move from fully skewed to symmetric models. 
	\item
	We have computed this convergence result for the particular case of the stable law, $h(\gamma)=\gamma^{-1/\alpha}$, but other conditionally Gaussian L\'{e}vy processes are likely to be amenable to similar analysis and we will study these in future work. Similarly, non-Gaussian versions may also be amenable to similar analysis provided convenient moment expressions are available for certain distributions. 
\end{enumerate}

\subsection{The L\'{e}vy-driven sde}
Recall that the required sde solution is given in Eq. (\ref{eq:sde})
and in particular we require a representation of the integral:
\[{\bf I}({\bf f}_t)=\int_{0}^{1}{\bf f}_tdW(u),\,\,\,\,\,\,
{\bf f}_t=e^{A(t-u)}{\cal I}(u\leq t)\mathbf{h}.\]
The proposed series representation of the integral is based on the generalised shot noise representation of the previous section,
\begin{equation}
W(t)=\sum_{i=1}^\infty h(\Gamma_i)U_i{\cal I}(V_i\leq t)-tc_i\label{Levy_series1}
\end{equation}
and to get the basic result we may (informally) substitute this directly into (\ref{eq:sde}) to obtain
\begin{IEEEeqnarray}{rCl}
{\bf I}({\bf f}_t)   {=}  
\sum_{i=1}^{\infty}\left[ \
U_i\Gamma_i^{-1/\alpha}{\bf f}_t(V_i)
-b_i^{(\alpha)}{E}[U_1]
{E}[{\bf f}_t(V_1)]
\right]
\IEEEeqnarraynumspace
\label{eq:int_approx}
\end{IEEEeqnarray}
where as before:
\begin{itemize} 
	\item $\{\Gamma_i \}$ are event times of a  unit-rate Poisson process 
	\item $\{V_i\}$ are i.i.d. $\mathcal{U}(0,1)$, 
	\item $\{U_i\}$ are i.i.d. such that ${E}[|U_i|^\alpha]<\infty$
	\item $\{b_i^{(\alpha)}\}$ are constants, non-zero only if $\alpha \in (1,2)$ 
\end{itemize}
and ${E}[{\bf f}_t(V_1)]$ is as follows:
\begin{align*}
{E}[{\bf f}_t(V_1)]&=\int_{0}^t e^{A(t-u)}\mathbf{h}du
\end{align*}
A similar form has been proved to converge to the correct $\alpha$-stable sde in \cite{Samorodnitsky_Taqqu_1994}, in the scalar sde case,  using a discrete Rademacher random variable $U_i\in\{-1,+1\}$. In \cite{Lemke_Godsill_2015} the scalar case was also proven with  $U_i\sim {\cal N}(\mu_W,\sigma_W^2)$, and a special case was also given for symmetric Gaussians in \cite{Samorodnitsky_Taqqu_1994} (the so-called Condition S). Here we are extending to the vector sde $X(t)\in \Re^P$, with  $U_i\sim {\cal N}(\mu_W,\sigma_W^2)$, for which the full proof will be presented in a future publication. 
Some insights may be gained by the following informal interpretation (see also \cite{Samorodnitsky_Taqqu_1994}):
\begin{itemize}
	\item $\{V_i\}$ are jump arrival times
	\item $U_i\Gamma_i^{-1/\alpha}$ is the size of jump at $V_i$
	\item ${\bf f}(V_i)$ is the effect at time $t$ of passing a unit jump at $V_i$ through the linear ODE $dX_t=AX_tdt$
	\item $b_i^{(\alpha)}{E}[W_1]
	{E}[{\bf f}(V_1)]$ is a drift term, only non-zero for $1< \alpha <2$
\end{itemize}
In particular, for practical implementations we here introduce the randomly 
truncated version, with truncation for $\Gamma_i\leq c$,
\[{\bf I}_c({\bf f}_t) = 
\left\{\sum_{i:\Gamma_i\leq c}
W_i\Gamma_i^{-1/\alpha}{\bf f}_t(V_i)\right\} -{\cal I}(\alpha>1)A(c){E}[{\bf f}_t(V_1)]
\]
which has zero mean for finite $c$ and $\alpha>1$, by construction. 

We now analyse the error due to this truncation, given by $R^c_t={\bf I}({\bf f}_t)-{\bf I}_c({\bf f}_t)$. The following result gives an exact characterisation of the covariance of the error:
\begin{thm}
	\begin{align*}\Sigma_t^c=\text{cov}(R^c_t)=(\sigma_W^2+\mu_W^2)\frac{\alpha}{2-\alpha}c^{1-2/\alpha}E[{\bf f}_t(V_1){\bf f}_t(V_1)^T].
	\end{align*}
\end{thm}
{\em Proof:\/}
The proof relies on calculation of the covariance of the point process corresponding to $R^c_t$ and will be presented in a future publication.

{\em Remark:\/}
The expectation is obtained for our sde as:
\[
E[{\bf f}_t(V_1){\bf f}_t(V_1)^T]=\int_0^t\exp(A(t-u)){\bf h}{\bf h}^T\exp(A^T(t-u))du
\]
and we can recognise this as the covariance function of a linear Gaussian sde \cite{Oksendal2003}
\[
dZ(t)=AZ(t)dt+{\bf h}dB(t)
\]  
where 
$\{B(t)\}$ is the unidimensional
Brownian motion.
Hence the covariance of $\{R^c_t\}$ is exactly matched by that of a linear Gaussian sde $\hat{R}_t^c$:
\[d\hat{R}^{c}_t=A\hat{R}^{c}_tdt+\left((\sigma_W^2+\mu_W^2)\frac{\alpha}{2-\alpha}c^{1-2/\alpha}\right)^{0.5}{\bf h}dB(t)\]
and $\hat{R}^{c}_0=0$. 

We propose then to approximate the exact integral ${\bf I}({\bf f}_t)$ in four possible ways, each involving a different truncation/ Gaussian approximation to the residual, as explained in the next section. 

\section{L\'{e}vy state space model}
We are now in a position to specify the stable L\'{e}vy state-space model. This requires the specification of the transition density between times $s$ and $t$, $t>s$, expressed as
\begin{equation}
{X}(t)=e^{A(t-s)}{X}(s)+\int_{s}^{t}e^{A(t-u)}\mathbf{h}dW(u).\label{eq:sde_soln}
\end{equation}
This may be obtained fairly straightforwardly from (\ref{eq:int_approx}) and the properties of the matrix exponential and the point process on the sub-interval $(s,t]$ of $[0,1]$ as
\begin{align}
X_t&=e^{A(t-s)}X_s\nonumber\\&+ \delta_t^{1/{\alpha}}\sum_{i=1}^{\infty}\left[ \
U_i\Gamma_i^{-1/\alpha}{\bf f}_t(V_i)
-b_i^{(\alpha)}{E}[U_1]
\int_{s}^t{\bf f}_t(u)du
\right]\label{eq:ss_untruncated}
\end{align}
where $\delta_t=t-s$ and
all terms are defined as before, except that now $V_i\overset{iid}{\sim}{\cal U}(s,t)$ are uniforms on the sub-interval $(s,t]$ and $\{\Gamma_i,U_i\} 
$ are drawn independently of the point process in any non-overlapping time intervals. Thus, for each successive time interval a new set of $\Gamma,U,V$s needs to be independently generated afresh in order to forward simulate the process for times $>t$. We have constructed the model in this way so that a causal forward simulation may be carried out without requiring a prior simulation of all the $\Gamma,U,V$ for the entire time axis as in (\ref{eq:int_approx}) .   

The randomly truncated version $X^c_{t}$ can be obtained similarly,
\begin{align}
X^c_{t}&=e^{A(t-s)}X_s+ Z^c_{(s,t]} -\bar{Z}^c_{(s,t]}\label{eq:trans_trunc}
\end{align}
with 
\[
Z^c_{(s,t]}=\delta_t^{1/{\alpha}}\sum_{i:\Gamma_i\leq c\delta_t}
U_i\Gamma_i^{-1/\alpha}{\bf f}_t(V_i)\]
and\[\bar{Z}^c_{(s,t]}={\cal I}(\alpha>1)\mu_W\frac{\alpha}{\alpha-1}c^{1-1/\alpha}\int_{s}^{t}{\bf f}_t(u)du.
\]
We can here interpret truncation parameter $c$ as the expected number of jumps per unit time. 
Note that, conditional upon $\{\Gamma_i,\,V_i\}$ the term $Z^c_{(s,t]}$ is Gaussian, 
\[
Z^c_{(s,t]}\sim{\cal N}(\mu_Wm^c_{(s,t]},\sigma_W^2 S^c_{(s,t]})
\]
with 
\begin{equation}{m}^c_{(s,t]}=\delta_t^{1/\alpha}\sum_{i:\Gamma_i\leq c\delta_t} \Gamma_i^{-1/\alpha}\exp(A(t-V_i)){\bf h}\label{eq:m_sum}\end{equation}
and 
\begin{equation}
{S}^c_{(s,t]}=\delta_t^{2/\alpha}\sum_{i:\Gamma_i\leq c\delta_t} \Gamma_i^{-2/\alpha}\exp(A(t-V_i)){\bf h}{\bf h}^T\exp(A(t-V_i))^T. \label{eq:S_sum}\
\end{equation}

To complete the state space model, recall that the residual terms of the series from $(c,\infty)$, ${R}^c_{(s,t]}$, are zero mean with covariance matrix:
\[
(\sigma_W^2+\mu_W^2)\Sigma^c_{\delta_t}
\]
leading to a decomposition of the exact state space model as:
\[
{X}_{t}=\exp(A\delta_t)X_s + Z^c_{(s,t]} -\bar{Z}^c_{(s,t]}+{R}^c_{(s,t]}.
\] 
We can now propose several approximations to this exact representation by either removing or approximating the residual term ${R}^c_{(s,t]}$. In each case ${R}^c_{(s,t]}$ is approximated with an independent  zero mean Gaussian term having covariance $\nu(\sigma_W,\mu_W)\Sigma^c_{\delta_t}$, where $\nu(\sigma_W,\mu_W)$ is a constant leading to different approximations, each corresponding to a  conditionally Gaussian form for the approximated L\'{e}vy-driven sde's transition density:\begin{align}
\hat{X}_{t}|\hat{X}_s,\{\Gamma_i,\,V_i&\}^c_{(s,t]}\nonumber\\\sim&{\cal N}(\exp(A\delta_t)\hat{X}_s + \mu_W{m}^c_{(s,t]} -\bar{Z}^c_{(s,t]},\nonumber\\&\,\,\,\,\,\,\,\,\,\,\,\,\,\,\,\,\,\,\,\,\,\,\sigma_W^2{S}^c_{(s,t]}+\nu(\sigma_W,\mu_W)\Sigma^c_{\delta_t}).
\end{align}
In each case the conditioning random variables are defined such that $\{\Gamma_i,\,V_i\}^c_{(s,t]}$ are generated {\em independently of all other time intervals that do not overlap with $(s,t]$\/} as:
\[
\Gamma_i=\Gamma_{i-1}+E_i,\,\,E_i\overset{iid}{\sim} {\cal E}(1),\,\Gamma_i\leq c\delta_t,\,\,V_i\overset{iid}{\sim} {\cal U}((s,t])
\]
where $ {\cal E}(1)$ denotes the exponential distribution with unit mean.

This conditional Gaussian transition density is the building block for all of the forward simulation and inference methods described subsequently.
The  four cases of approximation considered are:
\begin{enumerate}
	\item  {\bf Truncated series ${X}^c_{t}$.}
	The first approximation neglects the residual ${R}^c_{(s,t]}$ altogether, hence leading to the truncated version of the process, ${X}^c_{t}$ (\ref{eq:trans_trunc}), and corresponding to $\nu(\sigma_W,\mu_W)=0$. 
	This approximation of the process has some computational benefits, as will be considered shortly, but is most likely the least accurate of our three approximations. 
	\item{\bf Gaussian residual approximation $\hat{X}^c_{t}$.}
	In this more sophisticated  approximation a zero mean Gaussian approximation to the residual is made, replacing ${R}^c_{(s,t]}$  a zero mean Gaussian random variable $\hat{R}^c_{(s,t]}$ with covariance matched to ${R}^c_{(s,t]}$. In this case $\nu(\sigma_W,\mu_W)=\mu_W^2+\sigma_W^2$.  
	This is most likely the most accurate representation of the process that we have, but it has some less desirable computational properties, as will be described shortly. Note however that for the important symmetric case $\mu_W=0$ the transition density is identical to the following partial Gaussian approximation, acquiring its computational advantages (principally the ability to marginalise $\sigma_W^2$). 
	\item{\bf Partial Gaussian residual approximation $\tilde{X}^c_{t}$.} 
	Here a halfway house may be considered in which $\hat{R}^c_{(s,t]}$ is replaced with the solution of a linear SDE whose covariance equals $\sigma_W^2\Sigma^c_{\delta_t}$ and hence $\nu(\sigma_W,\mu_W)=\sigma_W^2$. This version retains the computational benefits of ${X}^c_{t}$, but does not fully account for the covariance of ${R}^c_{(s,t]}$ (since $\sigma_W^2\leq \sigma_W^2+\mu_W^2$), and so is likely to sit in between 1) and 2) in terms of approximation accuracy. 
	\item{\bf Joint Gaussian residual approximation of $m$ and $S$.}
	A final, more complex version of the truncation attempts to approximate the resisual series of $m$ and $S$ {\em jointly\/}. This has the same computational benefits as $X^c$ and $\tilde{X}^c$, but is a little more complex to compute and is not exposed here, see \cite{Riabiz2019} for full details.
\end{enumerate}

\subsection{Matrix-vector form for L\'{e}vy state space Model}
We are now in a position to state the model in standard linear (conditionally upon $m^c_{(s,t]}$ and $S^c_{(s,t]}$) Gaussian state-space form. For the approximated process $\hat{X_t}$, which may be set equal to ${X^c_t}$, case~1) above with  $\nu(\mu_W,\sigma_W)=0$, $\hat{X}^c_t$, case~2) above with  $\nu(\mu_W,\sigma_W)=\mu_W^2+\sigma_W^2$, or $\tilde{X}^c_t$, case~3) above with  $\nu(\mu_W,\sigma_W)=\sigma_W^2$. In cases 1) and 3)  it may be written using an extended state vector $\alpha_t$:
\begin{equation}
\alpha_{t}=\begin{bmatrix} \hat{X}_{t}\\\mu_W\end{bmatrix},\,\,\,\,
\alpha_{t}={\bf A}\alpha_{s}+{\bf B}{\bf e}_{s,t},\,\,\, {\bf e}_{s,t}\sim {\cal N}(0,C_e)\label{eq:alpha_t}
\end{equation}
with
\begin{align*}
{\bf A}=\begin{bmatrix}e^{A\delta_t}&m^c_{(s,t]}-\bar{Y}^c_{(s,t]}\\
{\bf 0}^T&1\end{bmatrix},\,\,\,\,\,
{\bf B}=\begin{bmatrix}{\bf I}_{P\times P}\\{\bf 0}_{1\times P}\end{bmatrix}
\end{align*}
where $\bar{Y}^c_{(s,t]}=\bar{Z}^c_{(s,t]}/\mu_W$ (which does not depend upon $\mu_W$),
and where
$
C_e=\sigma_W^2{S}^c_{(s,t]}+\nu(\sigma_W,\mu_W)\Sigma^c_{\delta_t}
$. 
For case 2), owing to the more complex covariance structure of the Gaussian residual approximation $\hat{X}^c_t$, the mean $\mu_W$ may not be directly included as a state variable and we have
\begin{equation}
\hat{X}^c_{t}={\bf A}'\hat{X}^c_{s}+{\bf e}_{t},\,\,\,\,{\bf A}'=e^{A\delta_t}\label{eq:X_t}\end{equation} and
$
{\bf e}_{t}\sim N\left(m^c_{(s,t]}\mu_W-{\bar{Z}^c_{(s,t]}},\sigma_W^2S^c_{(s,t]}+(\sigma_W^2+\mu_W^2)\Sigma^c_{\delta_t}\right)
$.

\section{Computation of likelihoods and state-conditionals}
Equations (\ref{eq:alpha_t}) and (\ref{eq:X_t}) can now form the basis for forward simulation of data from the L\'{e}vy state space model, or for inference about its states and parameters. 
Likelihoods and state conditional distributions may be computed using standard Bayesian recursions, implemented by the Kalman Filter, all conditioned upon the latent variables $\{m,S\}$. This means that Monte Carlo likelihood-based and Bayesian inference procedures such as Monte Carlo EM (MCEM), Markov chain Monte Carlo (MCMC) and Sequential Monte Carlo (SMC) may routinely be implemented.

For a concrete example, take the case 1) or 3) approximated models and suppose we partially observe $\alpha_t$, at times $\{t_i\}_{i=1}^N$, $0\leq t_i\leq T$, arranged in ascending order,
\[
y_t={\bf H}\alpha_t+V_t,\,\,\,V_t\sim{\cal N}(0,\sigma_W^2\kappa_V),\,\,\kappa_V\geq 0\,.
\]
For example, if ${\bf H}=\begin{bmatrix}1&{\bf 0}_{1\times P}\end{bmatrix}$ and $\kappa_V={\bf 0}$ we are fully observing the first component of $\alpha_t$. 
Then, initialise $t_0=0$, $y_0=0$, $\alpha_{0}\sim {\cal N}(a_{0|0},C_{0|0})$ with $a_{0|0}=\begin{bmatrix}{\bf 0}_{P\times 1}\\\mu_{\mu_W}\end{bmatrix}$, $C_{0|0}=\sigma_W^2\begin{bmatrix}{\bf 0}_{P\times P}&{{\bf 0}_{P\times 1}}\\{\bf 0}_{1\times P}&\kappa_W\end{bmatrix}$, where $\kappa_W\geq 0$ is a constant. 
Note that the covariance of the noise $\sigma_W^2\kappa_V$ and the prior covariance of $\mu_W$, $\kappa_W\sigma^2_{W}$, are both scaled relative to $\sigma_W^2$. This is to ensure greatest analytical tractability of posterior densities computed using Kalman filtering. A flat prior on $\mu_W$ is obtained as $\kappa_W\rightarrow \infty$, and a fully observed case with no observation noise is obtained with $\kappa_V={ 0}$.  Other more general prior structures can of course be incorporated, but this will be at the cost of a convenient conjugate posterior density for $\sigma_W^2$.  

In order to obtain the closed form result, follow a scheme similar to that in, for example, \cite{Harvey1990}. In this redefine the dynamical noise as 
\[
\tilde{\bf e}_{s,t}\sim {\cal N}\left(0,\tilde{C}_e\right)\\,
\]
where ${\tilde C}_e=C_e/\sigma_W^2$, which does not depend on $\sigma_W^2$ by construction in cases 1) and 3). Then, 
the Kalman filter computes \cite{Harvey1990},  for each $i=1,...,N$,
\begin{align*}
p(\alpha_{t_i}|y_{t_{1:i-1}},\sigma_W^2)&={\cal N}(a_{t_i|t_{i-1}},\sigma_W^2{\tilde C}_{t_i|t_{i-1}})\\
p(\alpha_{t_i}|y_{t_{1:i}},\sigma_W^2)&={\cal N}(a_{t_i|t_i},\sigma_W^2{\tilde C}_{t_i|t_i})\\
p(y_{t_i}|y_{t_{1:i-1}},\sigma_W^2)&={\cal N}(\hat{y}_{t_i},\sigma_W^2{\bf F}_{t_i})
\end{align*} 
where $a_{t_i|t_i}$, ${\tilde C}_{t_i|t_i}$ etc. are the Kalman filter output variables under the definition of the modified dynamical noise distribution $\tilde{\bf e}$, and $\hat{y}_{t_i}={\bf H}a_{t_i|t_{i-1}}$, ${\bf F}_{t_i}={\bf H}{\tilde C}_{t_i|t_{i-1}}{\bf H}^T+\kappa_V$.
From this the marginal likelihood is obtained as 
\begin{align*}
p(&y_{t_{1:N}}|\sigma_W^2)=\prod_{i=1}^N p(y_{t_i}|y_{t_{1:i-1}},\sigma_W^2)\\
&=-\frac{MN}{2}\log(2\pi)-\frac{N\log(\sigma_W^2)}{2}-\frac{1}{2}\sum_{i=1}^N\log|{\bf F}_{t_i}|-\frac{1}{2\sigma_W^2}E_N\end{align*}
where
\[
E_N=\sum_{i=1}^N w_{t_i}^T {\bf F}_{t_i}^{-1} w_{t_i},\,\,\,
w_{t_i}=y_{t_i}-\hat{y}_{t_i}\nonumber\]
and where $M$ is the dimension of the observation vector $y_{t}$.

The conjugate prior  for this form of likelihood is the inverted gamma distribution $p(\sigma_W^2)={\cal IG}(\alpha_W,\beta_W)$. Completing the (standard) conjugate analysis \cite{Bernardo_Smith_1994}, we have: 
\begin{align}
\log &p(y_{t_{1:N}})=-\frac{MN}{2}\log(2\pi)-\frac{1}{2}\sum_{i=1}^N\log|{\bf F}_{t_i}|\nonumber\\&\hspace*{0.5cm}+\alpha_W\log\beta_W-(\alpha_W+N/2)\log(\beta_W+E_N/2)\nonumber\\&\hspace*{1cm}+\log\Gamma(N/2+\alpha_W)-\log\Gamma(\alpha_W)\label{eq:joint_likelihood}\end{align}
\begin{align*}
p(\sigma_W^2|y_{t_{1:N}})&={\cal IG}(\alpha_W+N/2,\beta_W+{E_N}/{2})\\
p(\alpha_{t_N}|\sigma_W^2,y_{t_{1:N}})&={\cal N}(a_{t_N|t_N},\sigma_W^2{\tilde C}_{t_N|t_N})
\end{align*}
and we have, remarkably, all the tools to carry out marginal anaylsis of the data conditioned on the auxiliary variables $(m,S)$, and also {\em joint\/} conditional analysis for $(\sigma_W^2,\mu_W,\{X\})$. We stress that this full analysis is only available with the chosen prior structures, in particular the normal-inverse gamma conjugate priors for $(\sigma_W^2,\mu_W)$ and the special scaled form of the noise covariance matrices $\sigma_W^2\kappa_V$ and $\kappa_W\sigma^2_{W}$. Without these structures we will still be able to use the standard Kalman filter model to marginalise and infer $X$, but we will lose the closed form marginal expressions for $\sigma_W^2$ and/or $\mu_W$. 
A minor modification applies for the joint Gaussian residual approximation, case 4), in which once again the full Kalman filter calculations can be implemented with a modified noise covariance (not detailed here). Finally, for the Gaussian residual approximation $\hat{X}^c_t$, Case 2), the conjugate likelihood structure is lost owing to the term $\sigma_W^2+\mu_W^2$ in the transition density's noise covariance $C_e$. The standard Kalman filter can be used to infer/marginalise $X$, but not $\mu_W$ or $\sigma_W^2$. 
\section{Example: Langevin model}
In order to give a concrete application case, we work through the analysis for a 2-dimensional ($P=2$) Langevin model with similar structure to that used in \cite{Christensen2012},
\[
d\dot{x}(t)=\theta\dot{x}(t)dt+dW(t),\,\,\,\,\theta<0\] 
and work with a state vector:
\[
X(t)=[x(t) \,\, \dot{x}(t)]^T\,.
\]
The required system matrices for this model are 
\[
A=\begin{bmatrix}0&1\\0&\theta\end{bmatrix},\,\,\,{\bf h}=\begin{bmatrix}0\\1\end{bmatrix}.
\]
Then, 
\[
\exp(At)=\exp(\theta t)\begin{bmatrix}0&1/\theta\\0&1\end{bmatrix}+\begin{bmatrix}1&-1/\theta\\0&0\end{bmatrix}
\]
and 
\begin{align*}
\exp&(At){\bf h}{\bf h}^T\exp(At)^T=\exp(2\theta t)\begin{bmatrix}1/\theta^2&1/\theta\\1/\theta&1\end{bmatrix}\\&+\exp(\theta t)\begin{bmatrix}-2/\theta^2& -1/\theta\\-1/\theta &0\end{bmatrix}+\begin{bmatrix}1/\theta^2 &0\\0&0\end{bmatrix}
\end{align*}
which renders the calculation of  $m$ and $S$ (see (\ref{eq:m_sum}) and (\ref{eq:S_sum}) fairly straightforward:
\begin{align}{m}=\delta_t^{1/\alpha}\sum_{i=1}^\infty \Gamma_i^{-1/\alpha}\left(\exp(\theta (t-V_i))\begin{bmatrix}1/\theta\\1\end{bmatrix}+\begin{bmatrix}-1/\theta\\0\end{bmatrix}\right)\label{eq:m_sum_langevin}\end{align}
and 
\begin{align}
{S}&=\delta_t^{2/\alpha}\sum_{i=1}^\infty \Gamma_i^{-2/\alpha}\left(\exp(2\theta (t-V_i))\begin{bmatrix}1/\theta^2&1/\theta\\1/\theta&1\end{bmatrix}\right.\nonumber\\&+\left.\exp(\theta (t-V_i))\begin{bmatrix}2/\theta^2&-1/\theta\\-1/\theta&0\end{bmatrix}+\begin{bmatrix}1/\theta^2&0\\0&0\end{bmatrix}\right)\label{eq:S_sum_langevin}.
\end{align}
Also:
\begin{align*}
{E}[{\bf f}_t(V_1)]&=\frac{1}{T}\int_{0}^t e^{A(t-u)}\mathbf{h}du\\
&=\frac{1}{T}\int_{0}^t \left(\exp(\theta (t-u))\begin{bmatrix}1/\theta\\1\end{bmatrix}+\begin{bmatrix}-1/\theta\\0\end{bmatrix}\right)du\\
&=\frac{1}{T} \left(\frac{1}{\theta}(\exp(\theta t)-1)\begin{bmatrix}1/\theta\\1\end{bmatrix}+t\begin{bmatrix}-1/\theta\\0\end{bmatrix}\right)
\end{align*}
and:
\begin{align*}
&\int_0^t\exp(A(t-u){\bf h}{\bf h}^T\exp(A(t-u)^Tdu\\&=\frac{(\exp(2\theta t)-1)}{2\theta}\begin{bmatrix}1/\theta^2&1/\theta\\1/\theta&1\end{bmatrix}\\&+\frac{\exp(\theta t)-1}{\theta}\begin{bmatrix}-2/\theta^2& -1/\theta\\-1/\theta &0\end{bmatrix}+t\begin{bmatrix}1/\theta^2 &0\\0&0\end{bmatrix}.
\end{align*}
\subsection{Forward simulation example}
\begin{figure}[h]
	\centering
	\includegraphics[width=0.45\textwidth]{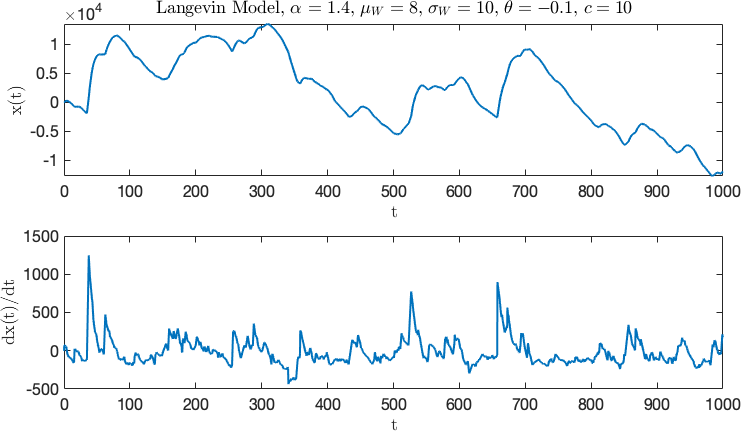}
	\caption{Example Langevin model data set, $\alpha=1.4$, positively skewed (\small{$\mu_W>0
		$)}}
	\label{alpha_1_4_data}
\end{figure}

\begin{figure}[h]
	\centering
	\includegraphics[width=0.45\textwidth]{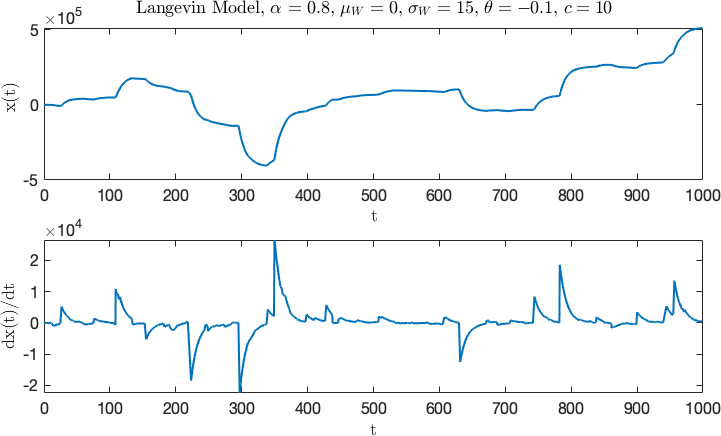}
	\caption{Example Langevin model data set, $\alpha=0.8$, symmetric ($\mu_W=0$)}
	\label{alpha_0_8_data}
\end{figure}

Equation (\ref{eq:alpha_t}) may be used directly to simulate a discrete `skeleton' of the process, and this in itself may be of use in the study of these intractable processes. An example simulation is given in Fig. \ref{alpha_1_4_data} for an asymmetric Langevin model with $\alpha=1.4$. 
A much more heavy-tailed example is given in Fig. \ref{alpha_0_8_data} for the symmetric case.


\subsection{Marginal Monte Carlo Filter}
\begin{figure}[h]
	\centering
	\includegraphics[width=0.45\textwidth]{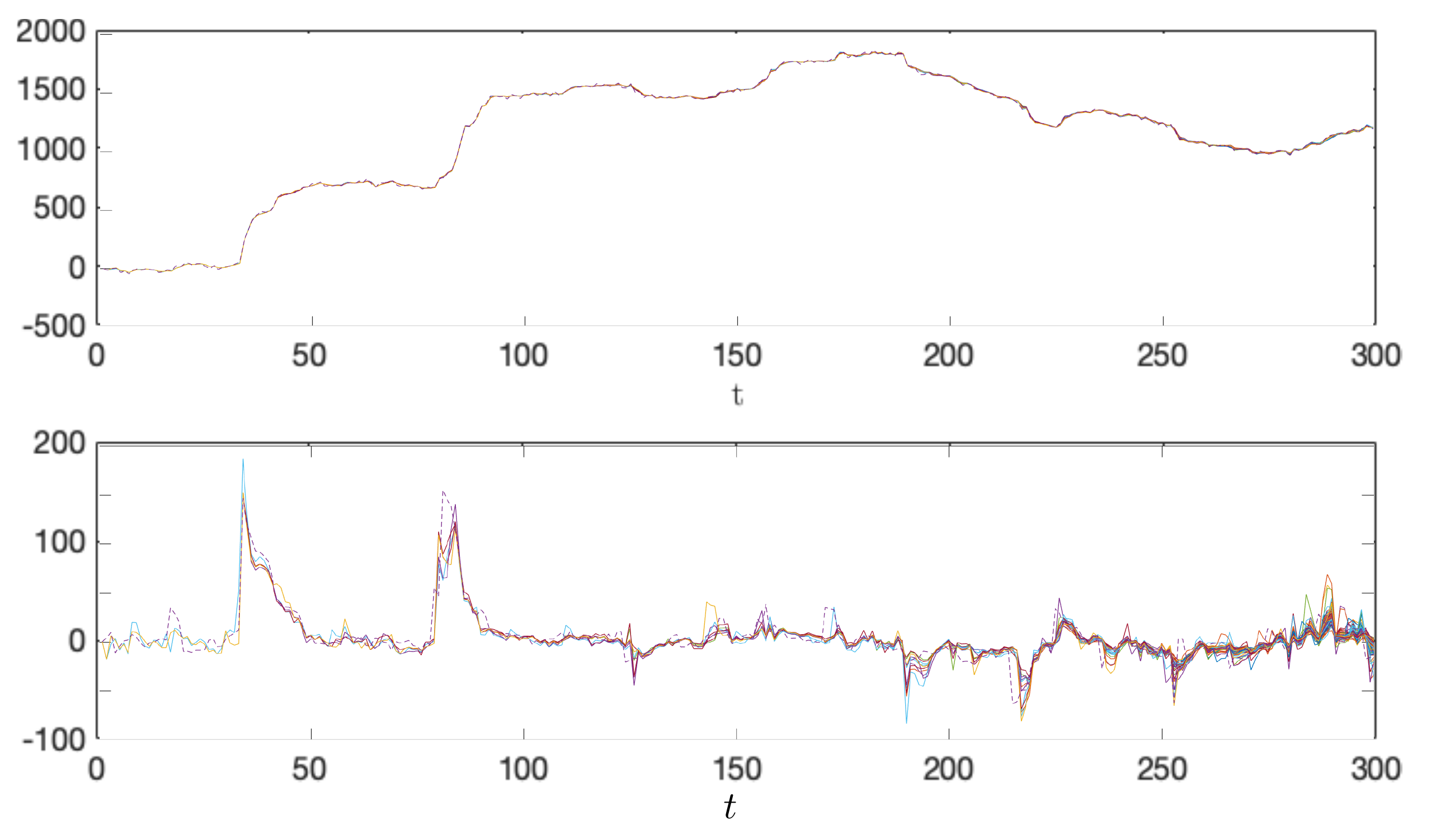}
	\caption{Marginal SMC Filter, $\alpha=1.4$. Top panel: MC samples of the trajectory of $x(t)$, observed data - dashed line. Lower panel: MC samples of the trajectory of $\dot{x}(t)$ with true state shown dashed. Same parameters as Fig.~\ref{alpha_1_4_data}}
	\label{alpha_1_4_data_MC}
\end{figure}

\begin{figure}[h]
	\centering
	\includegraphics[width=0.45\textwidth]{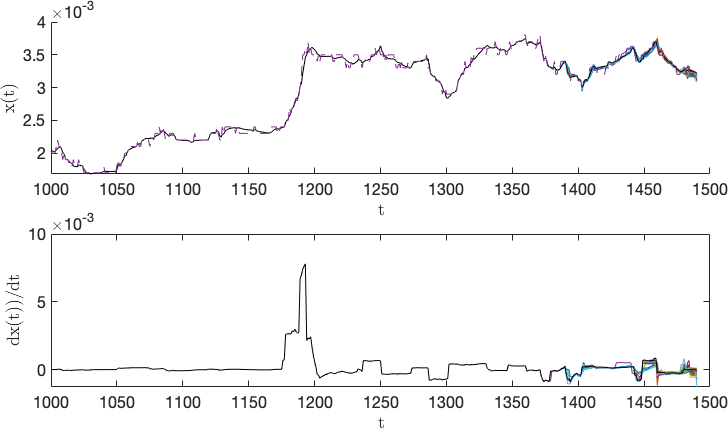}
	\caption{Processing of EURUSD exchange rate, `tick' data from 2006, $\alpha=0.8$, 4000 particles. Observed data shown dashed in top panel. Marginal SMC sample paths}
	\label{alpha_0_8_EURUSD}
\end{figure}

A sequential Monte Carlo (SMC) method is applied to infer the posterior distribution of the states of the Langevin model.
This is implemented in the form of a standard bootstrap particle filter, proposing the $m,S$ variables according to  (\ref{eq:m_sum_langevin}) and (\ref{eq:S_sum_langevin}) at each time interval and computing marginal weights according to (\ref{eq:joint_likelihood}). This is a fairly standard implementation, see e.g. \cite{Cappe2007}, and we omit further details here, leaving this to a subsequent publication. In all cases the state is observed as $y(t)=x(t)+V(t)$ and the derivative $\dot{x}(t)$ is unobserved. An example run of the filter, showing the Monte Carlo filter trajectories, is displayed in Fig. \ref{alpha_1_4_data_MC}. Note that $\dot{x}(t)$ is quite well inferred in the bottom panel, with appropriate posterior uncertainty about the trajectory. In Fig. \ref{alpha_0_8_EURUSD} a small segment of the high-frequency tick level Euro
dollar exchange rate from 2006 is modelled and tracked over a period of a few minutes with the filter. Here we found that a very heavy-tailed model ($\alpha=0.8$) was successfully able to capture the apparent rapid jumps in the trend ($\dot{x}$) of the process. Full simulation results and details of implementation will be provided in subsequent publications. 

\section{Conclusion}
Here we have presented a complete framework for modelling and inference in non-Gaussian L\'{e}vy-driven sdes. We have focussed on the $\alpha$- stable case, and we have outlined how the same methodological framework may be extended to more general classes of L\'{e}vy processes through different choices of the function $h()$, and more generally $H()$, in the shot noise representation. Further publications will present our proofs in full and explore further applications and other classes of L\'{e}vy process within our L\'{e}vy state space framework.


\bibliographystyle{unsrt}

{
\bibliography{Asilomar19.bib}
}

\end{document}